\documentclass[fleqn]{amsart}  % Use amsart for arxiv submission
\usepackage{graphicx}   % For including graphics
\usepackage{multirow}   % For multi-row tables
\usepackage{amsmath, amssymb, amsfonts}  % AMS packages for math
\usepackage{amsthm}     % For theorem environments
\usepackage{mathrsfs}   % For math script fonts
\usepackage[title]{appendix} % For appendix support
\usepackage{xcolor}     % For text color
\usepackage{textcomp}   % For text symbols
\usepackage{manyfoot}   % For multiple footnotes
\usepackage{booktabs}   % For nicer tables

\theoremstyle{plain} % Plain style for theorems
\newtheorem{theorem}{Theorem}[section] % Theorems numbered by section
\theoremstyle{definition} % Definition style

\newtheorem{remark}{Remark}[section]
\newtheorem{lem}{Lemma}[section]
\newtheorem{cor}{Corollary}[section]

\theoremstyle{definition} % Definition style
\newtheorem{dfn}{Definition}
\numberwithin{equation}{section} % Number equations by section

\begin{document}
	
	\title{Quaternionic Generalization of the Enestr\"{o}m-Kakeya Theorem}
	\author{D. Tripathi}\email{dinesh@mru.edu.in, dn3pathi.phdcpa@gmail.com}
	\address{Department of Sciences-Mathematics\\
	Manav Rachna University-121004, INDIA}
	
	\begin{abstract}
		{In 2020, Carney et.al. proved the quaternionic version of the Enestr\"{o}m-Kakeya Theorem, which states that a polynomial $p(q)=\sum_{\nu=0}^n q^\nu a_\nu$ with non-negative and monotonically increasing coefficients $(0<a_0\le a_1\le \cdots \le a_n)$ has all of its zeros within the unit ball $|q|\le 1$. Numerous generalizations of Enestr\"{o}m-Kakeya Theorem are available in the literatures (\cite{m}-\cite{mmr}). In this paper, we extend some of these generalizations to the quaternionic context and present several potential results.}
	\end{abstract}

	\maketitle
	

%{\bf Mathematics Subject Classification (2010):} 30B50, 30D30, 30E20.

\section{Introduction}\label{sec1}
Investigating polynomial zeros and their locations in the plane is a topic of interest in both mathematics and practical fields such as physical systems, utilizing various methods from geometric function theory. A substantial part of classical geometric function theory focuses on techniques to determine bounds for polynomial zeros. These methods are crucial for developing strategies in contemporary research. Much work has been devoted to identifying regions that contain all the zeros of a polynomial, which are often circular or annular in shape. Enestr\"{o}m and Kakeya \cite{m,mmr} were pioneers in this area, providing a classical solution for the polynomial with complex
variables and under constraints on the polynomial coefficients. Basically, they proved the following;
\\
\begin{theorem}\label{ek}
    Let $p(z)=\sum_{\nu=0}^{n}a_\nu z^\nu$ is a polynomial of degree $n$ with real coefficient satisfying $0< a_0\le a_1\le \cdots\le a_n$, then all the zeros of $p(z)$ lie in $|z|\le 1$.
\end{theorem}

This was followed by numerous studies exploring the zero distribution for polynomials with restricted coefficients, with comprehensive overviews available in the works of Marden \cite{m} and Milovanovi\'{c} {\it et al.}\cite{mmr}, Tripathi \cite{dt1}, etc.

Let $$\mathbb{H}:=\{q;~~ q=x_0+i x_1+j x_2+k x_3,~~\forall x_i\in \mathbb{R}, i=0, 1, 2, 3\},$$ where $i, j, k$ satisfy
$$i^2=j^2=k^2=-1, ij=-ji=k, jk=-kj=i, ki=-ik=j.$$
The set $\mathbb{H}$ be non-commutative ring of quaternions. By $q\in \mathbb{H}$, we mean the composition of real part $\operatorname{Re}(q)=x_0$ and $\operatorname{Im}(q)=x_1 i+x_2 j+x_3 k$. We denote the quaternionic conjugate of $q$ as $\Bar{q}:=x_0-x_1 i-x_2 j-x_3 k$ and the norm of $q$ is $|q|=\sqrt{q\bar{q}}=\sqrt{x_0^2+x_1^2+x_2^2+x_3^2}$. The inverse of $q\ne0 \in \mathbb{H}$ is given by $q^{-1}=|q|^{-2}q.$
For $r>0$, define the open ball $$\mathcal{B}(0, r)=\{q\in \mathbb{H}; |q|<r\}$$ and $\mathbb{B}=\{q; |q|<1\}$ denotes the open unite ball in $\mathbb{H}$. The set $$\mathbb{S}:=\{q=x_1 i+x_2j+x_3k;~ x_1^2+x_2^2+x_3^2=1\}$$ is the unit sphere of purely imaginary quaternions.

Let $q$ be the indeterminate for a quaternionic polynomial $$p(q)=\sum_{\nu=0}^n q^\nu a_\nu,~~ a_n\ne 0,~~ a_\nu \in \mathbb{H},~~ 0\le \nu\le n$$ of degree $n$. The polynomial $p(q)$ satisfy the regularity conditions and their behavior resembles very
closely to that of holomorphic functions of a complex variable.  A different product $(\star)$ of such kind of polynomials define over the skew-field $\mathbb{H}$ as
\bigskip

\begin{dfn}\cite[Definition 3.1]{gs2}
    If $f(q)=\sum_{\nu=0}^{\infty}q^\nu a_\nu$ and $g(q)=\sum_{\nu=0}^{\infty}q^\nu b_\nu$ are given quaternionic power series with radii of convergence greater than $R$. Then, the real product $f$ and $g$ is $f\star g(q)$ and defined as $$f\star g(q)=\sum_{\nu=0}^{\infty}q^\nu c_\nu,~\textnormal{where}~~ c_\nu=\sum_{\mu=0}^{\nu}a_\mu b_{\nu-\mu}.$$
\end{dfn}

Recently, an analytic theory of functions of a quaternionic variable has been established \cite{gs1, gs2}. Specifically, a maximum modulus theorem for regular functions of a quternionic variable has been introduced by Gentili and Struppa \cite{gs1}, which includes convergence power series and polynomials in the form of the following result.
\\
\begin{theorem}\label{mm}
    (Maximum Modulus Theorem). Let $\mathcal{F}:\mathcal{B}\to \mathbb{H}$ be a regular function, where $\mathcal{B} = \mathcal{B}(0, r)$ be a ball in $\mathbb{H}$ with center at origin and radius $r > 0$. If $|\mathcal{F}|$ has a relative maximum at a point $a\in \mathcal{B}$, then $\mathcal{F}$ is a constant on $\mathcal{B}$.
\end{theorem}

 In 2020, Carney {\it et al.} \cite{cb} extended Theorem \ref{ek} for the class of polynomials $p(q)$, which is as follows;
\\
 \begin{theorem}\label{q1}
     Let $p(q)=\sum_{\nu=0}^{n} q^\nu a_\nu$ be a polynomial of degree $n$ with real coefficient satisfying $0<a_0\le a_1\le \cdots\le a_n$, then all the zeros of $p(q)$ lie in $|q|\le 1$.
 \end{theorem}
  
In the same paper Carney {\it et al.} \cite{cb} proved the following extension of the results due to Govil {\it et. al} \cite{gr} for the polynomial $p(q)=\sum_{\nu=0}^n q^\nu a_\nu$.
\\

\begin{theorem}\label{q2}
    Let $p(q)=\sum_{\nu=0}^n q^\nu a_\nu$ be the polynomial of degree $n$. Let $b\ne 0$ be a non-zero quaternion and $\measuredangle(a_\nu, b)\le \theta \le \pi/2$ for some $\theta$ and $0\le \nu \le n$. Assume $$|a_0|\le |a_1|\le \cdots \le|a_n|.$$ Then all the zeros of $p(q)$ lie in $$|q|\le \cos{\theta}+\sin{\theta}+2\frac{\sin{\theta}}{|a_n|}\sum_{\nu=0}^{n-1}|a_\nu|.$$
\end{theorem}
\begin{theorem}\label{q3}
    If $p(q)=\sum_{\nu=0}^{n} q^\nu a_\nu$ is a polynomial of degree $n$ where $a_\nu=\alpha_\nu+\beta_\nu i+\gamma_\nu j+\delta_\nu k$ for $0\le \nu\le n,$ satisfying $$0\le \alpha_1 \le \alpha_2\le \cdots\le \alpha_n, ~~\alpha_n\ne 0.$$ Then all the zeros of $p(q)$ lie in $|q|\le 1+\frac{2}{\alpha_n}\sum_{\nu=0}^n |a_\nu|.$
\end{theorem}
\bigskip

Several extensive generalizations of Theorem \ref{q1}, Theorem \ref{q2} and Theorem \ref{q3}  have been developed in recent years, including the works by Garden and Taylor \cite{gt}, Tripathi \cite{dt} etc. The objective of this paper is to prove several generalizations of Theorem \ref{q1} for a different class of polynomials.

\section{Main Results}\label{sec2}
Let $0=n_0<n_1<n_2<\cdots<n_k=n$ and let
\begin{eqnarray}
    \mathbb{P}_n:=\left\{p~;~p(q)=\sum_{\nu=0}^{k} q^{n_\nu} a_{n_\nu}, a_{n_\nu}\in \mathbb{H}, \nu=0, 1, 2, \ldots, k\right\}
\end{eqnarray}
be a set of polynomial over $\mathbb{H}$. Then our first result of this paper stated as;
\\
\begin{theorem}\label{t1}
   Let $p\in \mathbb{P}_n$ with $0<|a_{n_0}|\le|a_{n_1}|\le \cdots |a_{n_k}|$ such that no two adjacent coefficients $a_{n_j}'s$ are equal, i.e.
        $a_{n_j}\ne a_{n_{j+1}}, j=0, 1, 2, \cdots, k$ and for certain real numbers $\alpha$ and $\beta$ $$|\arg a_{n_j}-\beta|\le \alpha\le \pi/2, ~ j=0, 1, 2,\cdots, k.$$ Further let 
    \begin{eqnarray}\label{te1}    M_{n_j}:=\max_{|q|=1}\left|\frac{q^{n_{j}}a_{n_{j}}-q^{n_{j-1}+1}a_{n_{j-1}}}{a_{n_j}-a_{n_{j-1}}}\right|
    \end{eqnarray}
    with $$M:=\max_{1\le j\le k} M_{n_j}.$$
    Then all the zeros of $p(q)$ lie in 
    $$|q|\le M\left\{|a_n|(\cos{\alpha}+\sin{\alpha})+2\sum_{j=0}^{k-1}|a_{n_j}|\sin{\alpha}\right\}.$$
\end{theorem}

Our next result is the generalizations of Theorem \ref{q3} for the class of polynomials $\mathbb{P}_n$.
\\
\begin{theorem}\label{t2}
   Let $p\in \mathbb{P}_n$ with $a_{n_l}=\alpha_{n_l}+\beta_{n_l}i+\gamma_{n_l}j+\delta_{n_l}k$ for $l=0, 1, 2, \ldots, k$ such that no two adjacent coefficients $a_{n_j}'s$ are equal, i.e.
        $a_{n_l}\ne a_{n_{l+1}}, l=0, 1, 2, \cdots, k$ and 
        
        $$\alpha_{n_0}\le \alpha_{n_1}\le \cdots\le \alpha_{n_k},~ \beta_{n_0}\ge \beta_{n_1}\ge \cdots \ge \beta_{n_k},$$  $$\gamma_{n_0}\le \gamma_{n_1}\le \gamma_{n_{r-1}}\le \cdots \le \gamma_{n_r}\ge \gamma_{n_{r+1}}\ge \cdots \ge \gamma_{n_k}, 0\le r\le k $$ Further let $M_{n_j}, 1\le j\le k $ with $M$ as defined in Theorem \ref{t1}.
    Then all the zeros lie in 
    \begin{eqnarray}
       |q|&\le& \frac{M}{|a_n|}\left\{ (\alpha_{n}-\alpha_{0})+(\beta_{0}-\beta_{n})+2\gamma_{n_r}-(\gamma_{0}+\gamma_{n})\right.\nonumber\\
    &&\left.+2\sum_{j=0}^{k-1}|\delta_{n_j}|+|\delta_{n}|+|\alpha_{0}|+|\beta_{0}|+|\gamma_{0}|\right\}.
    \end{eqnarray}
\end{theorem}

\section{Lemma}
The following lemmas are required for the proof of our results.
\\
\begin{lem}\label{l1}
    Let $p(q)=\sum_{\nu=0}^n q^\nu a_\nu$ be a polynomial of degree $n$. Then for $R>1$
    \begin{eqnarray}\label{le1}
        \max_{|q|=R}|p(q)|\le R^n \max_{|q|=1}|p(q)|
    \end{eqnarray}
    with equality holds for $p(q)=q^n\lambda$.
\end{lem}

\begin{proof}
    Let $\eta(q)$ be a polynomial of degree $n$, then for $|q|>\rho>0$, then for any $r>\rho$, we have
    \begin{eqnarray}\label{lp1}
        \max_{|q|\ge r}|\eta(q)|=\max_{|q|=r}|\eta(q)|.
    \end{eqnarray}
    Define $F(q)=\eta(q^{-1})$. Then by Maximum Modulus Theorem, we have for $|q|\le \frac{1}{r}$,
    \begin{eqnarray*}
        \max_{|q|\le 1/r}|F(q)|=\max_{|q|=1/r}|F(q)|,
    \end{eqnarray*}
    or equivalently obtained
    \begin{eqnarray}\label{lp2}
        \max_{|q^{-1}|\le 1/r}|\eta(q^{-1})|=\max_{|q^{-1}|=1/r}|\eta(q^{-1})|,
    \end{eqnarray}
    which is similar to (\ref{lp1}).
    Now for $R\ge 1$, we obtain
    \begin{eqnarray*}
        \max_{|q|=R}|\eta(q)|\le \max_{|q|=1}|\eta(q)|.
    \end{eqnarray*}
    By Considering $\eta(q)=q^{-n}\star p(q)$, we immediately obtained (\ref{le1}).
\end{proof}

Then next lemma is due to Carney {\it et.al} \cite{cb}.
\\
\begin{lem}\label{l2}
    Let $q_1, q_2\in \mathbb{H}$ where $\measuredangle(q_1, q_2)=2\theta'\le 2\theta$, and $|q_1|\le|q_2|$. Then
    \begin{eqnarray}
        |q_2-q_1|\le(|q_2|-|q_1|)\cos{\theta}+(|q_1|+|q_2|)\sin{\theta}.
    \end{eqnarray}
\end{lem}
The following result completely describes the zero sets of a regular product of two polynomials in terms of the zero sets of the two factors.
\\
\begin{lem}\cite{gs}\label{zr}
	Let $F$ and $G$ be given quaternionic polynomials. Then the convolution product $(F\star G)(q_0)=0$ if and only if $F(q_0) =0$ or $F(q_0)\ne0 \implies G(F(q_0)^{-1}q_0F(q_0))=0$.
\end{lem}

\section{Proof of Results:}
\begin{proof}[Proof of Theorem \ref{t1}]
    Define 
    \begin{eqnarray}\label{tp1}
        \xi(q)&=&p(q)\star(1-q)\nonumber\\
        &&=-q^{n+1}a_n+\sum_{j=1}^{k}(q^{n_{j}}a_{n_{j}}-q^{n_{j-1}+1}a_{n_{j-1}})+q^{n_{0}}a_{n_{0}}\nonumber\\
        &&=-q^{n+1}a_n+\xi_1(q),
    \end{eqnarray}
    where 
    \begin{eqnarray}\label{tp2}
        \xi_1(q)=\sum_{j=1}^{k}(q^{n_{j}}a_{n_{j}}-q^{n_{j-1}+1}a_{n_{j-1}})+q^{n_{0}}a_{n_{0}}.
    \end{eqnarray}
From Lemma \ref{zr}, $\xi(q)=p(q)\star (1-q)=0$ if and only if either $p(q)=0$ of $p(q)\ne 0\implies p(q)^{-1}qp(q)-1=0$, i.e., $p(q)^{-1}qp(q)=1.$ If $p(q)\ne 1$, then $q=1$. Therefore, the only zeros of $\xi(q)$ are $q=1$ and the zeros of $p(q)$. \\

Now for $|q|=R(>1)$, we get from (\ref{tp2})
\begin{eqnarray}\label{tp3}
    |\xi_1(q)|&\le& \sum_{j=1}^{k}\left|\frac{q^{n_{j}}a_{n_{j}}-q^{n_{j-1}+1}a_{n_{j-1}}}{a_{n_j}-a_{n_{j-1}}}\right||a_{n_j}-a_{n_{j-1}}|+|a_{n_0}|R^{n_0}
    \end{eqnarray}
    using Lemma \ref{l1} for $|q|=R(>1)$
\begin{eqnarray}\label{tp4}
    |\xi_1(q)|&\le& \sum_{j=1}^{k} M_{n_j} R^{n_j}|a_{n_j}-a_{n_{j-1}}|+|a_{n_0}|R^{n_0},
\end{eqnarray}
where $M_{n_j}$ is define in (\ref{te1}).
From (\ref{te1}) $$M_{n_j}\ge 1,~~ 1\le j\le k,$$ therefore $$\max_{1\le j\le k} M_{n_j}=M\ge 1.$$
So from (\ref{tp4}), we obtain
\begin{eqnarray}\label{tp5}
   |\xi_1(q)| \le M R^n \sum_{j=1}^{k}|a_{n_j}-a_{n_{j-1}}|+|a_{n_0}|R^{n_0}.
\end{eqnarray}
Using Lemma \ref{l2}, we get for $|q|=R(>1)$
\begin{eqnarray}\label{tp6}
    |\xi_1(q)|&\le& M R^n \left\{\sum_{j=1}^{k}(|a_{n_j}|-|a_{n_{j-1}}|)\cos {\alpha}+\sum_{j=1}^{k}(|a_{n_{j-1}}|+|a_{n_j}|)\sin{\alpha}\right.\nonumber\\
    &&\left.+|a_{n_0}|(\cos{\alpha}+\sin{\alpha})\right\}\nonumber\\
    %&\le& MR^n\left\{(|a_{n_k}|-|a_{n_0}|)\cos{\alpha}+(|a_{n_0}|+2\sum_{j=1}^{k-1}|a_{n_j}|+|a_{n_k}|)\sin{\alpha}+|a_{n_0}|\right\}\nonumber\\
    &\le&M R^n\left\{|a_n|(\cos{\alpha}+\sin{\alpha})+2\sum_{j=0}^{k-1}|a_{n_j}|\sin{\alpha}\right\}
\end{eqnarray}
From (\ref{tp1}), we get for $|q|=R(>1)$
\begin{eqnarray}\label{tp7}
    |\xi(q)|&\ge& |a_n|R^{n+1}-|\xi_1(q)|\nonumber\\
    &\ge&|a_n|R^{n+1}-M R^n\left\{|a_n|(\cos{\alpha}+\sin{\alpha})+2\sum_{j=0}^{k-1}|a_{n_j}|\sin{\alpha}\right\}\nonumber\\
    &>& 0
\end{eqnarray}
for
$$R>M\left\{|a_n|(\cos{\alpha}+\sin{\alpha})+2\sum_{j=0}^{k-1}|a_{n_j}|\sin{\alpha}\right\}.$$
Which follows Theorem \ref{t1}.
\end{proof}
\begin{remark}
    We also have the inequality (\ref{tp5}) by removing the restriction of the equality of two adjacent coefficients from some $a_{n_j}, 0\le j\le n$ of Theorem \ref{t1}, i.e. we have inequality (\ref{tp5}) by taking $$a_{n_0}=a_{n_1}\ne a_{n_2}\ne a_{n_3}=a_{n_4}=a_{n_5}\ne a_{n_6}\ne a_{n_7}\ne \cdots \ne a_{n_{k-2}}=a_{n_{k-1}}\ne a_{n_k}.$$
    For that, we will define
    \begin{eqnarray}\label{re1}
        M_{n_1}&=&\max_{|q|=1}\left|\frac{a_{n_0}-q^{n_0+1}a_{n_0}+q^{n_1}a_{n_1}}{a_{n_0}-a_{n_0}+a_{n_1}}\right|\ge 1\\
        M_{n_2},&&~~~~~~~ {\textnormal{(same as (\ref{te1}) and not consider $M_{n_3}, M_{n_4}$ )}}\\
    \end{eqnarray}
    \begin{eqnarray}
      M_{n_5}&=&\max_{|q|=1}\left|\frac{q^{n_5}a_{n_5}-q^{n_4+1}a_{n_4}+q^{n_4}a_{n_4}-q^{n_3+1}a_{n_3}+q^{n_3}a_{n_3}-q^{n_2+1}a_{n_2}}{a_{n_5}-a_{n_4}+a_{n_4}-a_{n_3}+a_{n_3}-a_{n_2}}\right|~~~~~~\nonumber\\
        &\ge& 1,\\
        && M_{n_6}, M_{n_7},\cdots, M_{n_{k-2}}~~\textnormal{(same as (\ref{te1}) and not consider$M_{n_{k-2}}$)}
    \end{eqnarray}
        
        \begin{eqnarray}
        M_{n_{k-1}}&=&\max_{|q|=1}\left|\frac{q^{n_{k-1}}a_{n_{k-1}}-q^{n_{k-2}+1}a_{n_{k-2}}+q^{n_{k-2}}a_{n_{k-2}}-q^{n_{k-3}+1}a_{n_{k-3}}}{a_{n_{k-1}}-a_{n_{k-2}}+a_{n_{k-2}}-a_{n_{k-3}}}\right|\nonumber\\
        &\ge&1\\
        &&M_{n_k} \textnormal{(same as (\ref{te1}))}\\~ \label{re2}
        M&=&\max(M_{n_1}, M_{n_2}, M_{n_5}, M_{n_6},\cdots, M_{n_{k}})\nonumber\\
       && \textnormal{(instead of the expression $M$ as defined in Theorem \ref{t1})}~~~~~~~\nonumber\\
        &\ge&1
    \end{eqnarray}
    Then form (\ref{tp1})
    \begin{eqnarray}
        \xi(q)&=&-q^{n+1}a_n+(q^{n_{k}}a_{n_{k}}-q^{n_{k-1}+1}a_{n_{k-1}})\nonumber\\
        &&+(q^{n_{k-1}}a_{n_{k-1}}-q^{n_{k-2}+1}a_{n_{k-2}}+q^{n_{k-2}}a_{n_{k-2}}-q^{n_{k-2}+1}a_{n_{k-2}})\nonumber\\
        &&+(q^{n_{k-3}}a_{n_{k-3}}-q^{n_{k-4}+1}a_{n_{k-4}})+\cdots+(q^{n_6}a_{n_6}-q^{n_5+1}a_{n_5})\nonumber\\
        &&+(q^{n_5}a_{n_5}-q^{n_4+1}a_{n_4}+q^{n_4}a_{n_4}-q^{n_3+1}a_{n_3})+(q^{n_3}a_{n_3}-q^{n_2+1}a_{n_2})\nonumber\\
        &&+(q^{n_2}a_{n_2}-q^{n_1+1}a_{n_1})+q^{n_{0}}a_{n_{0}}\nonumber\\
        &=&-q^{n+1}a_n+\xi_1(q).
    \end{eqnarray}
    On taking account of the values of $M_{n_j}$ from (\ref{re1})-(\ref{re2}) and using it as same as the proof of Theorem \ref{t1} to obtain for $|q|=R(>1)$,
    \begin{eqnarray}
        |\xi_1(q)|&\le& |a_{n_k}-a_{n_{k-1}}|M_{n_k}R^{n_k}+|a_{n_{k-1}}-a_{n_{k-2}}+a_{n_{k-2}}-a_{n_{k-3}}|M_{n_{k-1}}R^{n_{k-1}}\nonumber\\
        &&+\sum_{j=6}^{k-3}|a_{n_j}-a_{n_{j-1}}|M_{n_j}R^{n_j}+|a_{n_5}-a_{n_4}+a_{n_4}-a_{n_3}|M_{n_5}R^{n_5}\nonumber\\
        &&+|a_{n_3}-a_{n_2}|M_{n_3}R^{n_3}+|a_{n_2}-a_{n_1}|M_{n_2}R^{n_2}+|a_{n_0}|R^{n_0}\nonumber\\
        &\le& MR^{n}\left\{|a_{n_k}-a_{n_{k-1}}|+|a_{n_{k-1}}-a_{n_{k-2}}+a_{n_{k-2}}-a_{n_{k-3}}|\right.\nonumber\\
        &&\left.+\sum_{j=6}^{k-3}|a_{n_j}-a_{n_{j-1}}|+|a_{n_5}-a_{n_4}+a_{n_4}-a_{n_3}|+|a_{n_3}-a_{n_2}|\right.\nonumber\\
        &&\left.+|a_{n_2}-a_{n_1}|\right\}+|a_{n_0}|R^{n_0}\nonumber\\
        &\le&MR^n\left\{\sum_{j=2}^k|a_{n_j}-a_{n_{j-1}}|+|a_{n_0}|\right\}.
    \end{eqnarray}
    
\end{remark}

\begin{proof}[Proof of Theorem \ref{t2}]
    From inequality (\ref{tp31}) of the proof of Theorem \ref{t1}, we have for $|q|=R(>1)$
\begin{eqnarray}\label{tp21}
    |\xi_1(q)|&\le& M R^n \left\{ \sum_{j=1}^{k}|a_{n_j}-a_{n_{j-1}}|+|a_{n_0}|\right\}\nonumber\\
    &=& M R^n \left\{ \sum_{j=1}^{k}|(\alpha_{n_j}-\alpha_{n_{j-1}})+i(\beta_{n_j}-\beta_{n_{j-1}})+j(\gamma_{n_j}-\gamma_{n_{j-1}})\right.\nonumber\\
    &&\left.+k(\delta_{n_j}-\delta_{n_{j-1}})|+|a_{n_0}|\right\}\nonumber\\
    &\le&M R^n \left\{ \sum_{j=1}^{k}(|\alpha_{n_j}-\alpha_{n_{j-1}}|+|\beta_{n_j}-\beta_{n_{j-1}}|+|\gamma_{n_j}-\gamma_{n_r}+\gamma_{n_r}- \gamma_{n_{j-1}}|\right.\nonumber\\
    &&\left.+|\delta_{n_j}-\delta_{n_{j-1}}|)+|\alpha_{n_0}|+|\beta_{n_0}|+|\gamma_{n_0}|+|\delta_{n_0}|\right\}
    \end{eqnarray}
    \begin{eqnarray}\label{tp22}
    &&\le M R^n \left\{ (\alpha_{n_k}-\alpha_{n_0})+(\beta_{n_0}-\beta_{n_k})+\sum_{j=1}^{k}(|\gamma_{n_j}-\gamma_{n_r}|+|\gamma_{n_r}- \gamma_{n_{j-1}}|)\right.\nonumber\\
    &&\left.+2\sum_{j=0}^{k-1}|\delta_{n_j}|+|\delta_{n_k}|+|\alpha_{n_0}|+|\beta_{n_0}|+|\gamma_{n_0}|\right\}\nonumber\\
    &&\le M R^n \left\{ (\alpha_{n_k}-\alpha_{n_0})+(\beta_{n_0}-\beta_{n_k})+2\gamma_{n_r}-(\gamma_{n_0}+\gamma_{n_k})\right.\nonumber\\
    &&\left.~~~~+2\sum_{j=0}^{k-1}|\delta_{n_j}|+|\delta_{n_k}|+|\alpha_{n_0}|+|\beta_{n_0}|+|\gamma_{n_0}|\right\}
\end{eqnarray}
Further from (\ref{tp1}), we get for $|q|=R$
\begin{eqnarray}\label{tp23}
    |\xi(q)|\ge R^{n+1}|a_n|-|\xi_1(q)|.
\end{eqnarray}
Using (\ref{tp22}), we have
\begin{eqnarray}
    |\xi(q)|&\ge&R^{n+1}|a_n|-M R^n \left\{ (\alpha_{n_k}-\alpha_{n_0})+(\beta_{n_0}-\beta_{n_k})+2\gamma_{n_r}\right.\nonumber\\
    &&\left.-(\gamma_{n_0}+\gamma_{n_k})+2\sum_{j=0}^{k-1}|\delta_{n_j}|+|\delta_{n_k}|+|\alpha_{n_0}|+|\beta_{n_0}|+|\gamma_{n_0}|\right\}\nonumber\\
    >0,
\end{eqnarray}
if 
\begin{eqnarray*}
    R&>&\frac{M}{|a_n|}\left\{ (\alpha_{n_k}-\alpha_{n_0})+(\beta_{n_0}-\beta_{n_k})+2\gamma_{n_r}-(\gamma_{n_0}+\gamma_{n_k})\right.\nonumber\\
    &&\left.+2\sum_{j=0}^{k-1}|\delta_{n_j}|+|\delta_{n_k}|+|\alpha_{n_0}|+|\beta_{n_0}|+|\gamma_{n_0}|\right\}
\end{eqnarray*}
and Theorem \ref{t2} follows.
\end{proof}
By taking $\gamma_{n_j}=0=\delta_{n_j}, 0\le j\le k$ in inequality (\ref{tp21}), we get for $|q|=R(>1)$
\begin{eqnarray*}
    |\xi_1(q)|&\le& M R^n \left\{ (\alpha_{n_k}-\alpha_{n_0})+(\beta_{n_0}-\beta_{n_k})+|\alpha_{n_0}|+|\beta_{n_0}|\right\}
\end{eqnarray*}
and hence obtain following generalization of Theorem \ref{q3}.\\
\begin{cor}\label{co1}
    Let $p\in \mathbb{P}_n$ with $a_{n_l}=\alpha_{n_l}+\beta_{n_l}i$ such that no two adjacent coefficients $a_{n_j}'s$ are equal, i.e.
        $a_{n_j}\ne a_{n_{j+1}}, j=0, 1, 2, \cdots, k$ and 
        $$\alpha_{n_0}\le \alpha_{n_1}\le \cdots\le \alpha_{n_k},~ \beta_{n_0}\ge \beta_{n_1}\ge \cdots \ge \beta_{n_k}.$$ Further let $M_{n_j}, 1\le j\le k $ with $M$ as defined in Theorem \ref{t1}.
    Then all the zeros lie in 
    \begin{eqnarray}
       |q|&\le& \frac{M}{|a_n|}\left\{ (\alpha_{n}-\alpha_{0}+|\alpha_{0}|)+(\beta_{0}-\beta_{n}+|\beta_{0}|)\right\}.
    \end{eqnarray}
\end{cor}

On removing the restriction from $\beta_{n_j}, j=0, 1, 2,\ldots,n$ in Corollary \ref{co1} i.e. $\alpha_{n_0}\le \alpha_{n_1}\le \cdots\le \alpha_{n_k}$ and for all $\beta_{n_j}$ in $p(q)=\sum_{\nu=0}^{k} q^{n_\nu} a_{n_\nu}$, we have from (\ref{tp21})
\begin{eqnarray}\label{r2}
    |\xi_1(q)|&\le& M R^n \left\{ \sum_{j=1}^{k}(|\alpha_{n_j}-\alpha_{n_{j-1}}|+|\beta_{n_j}-\beta_{n_{j-1}}|)+|\alpha_{n_0}|+|\beta_{n_0}|\right\}\nonumber\\
    &&\le M R^n \left\{ (\alpha_{n_k}-\alpha_{n_0}+|\alpha_{n_0}|)+2\sum_{j=0}^{k-1}|\beta_{n_j}|\right\}\nonumber\\
\end{eqnarray}
and on combining it with (\ref{tp23}), we have the following;
\\
\begin{cor}\label{co2}
   Let $p\in \mathbb{P}_n$ with $a_{n_l}=\alpha_{n_l}+\beta_{n_l}i $ such that no two adjacent coefficients $a_{n_j}'s$ are equal, i.e.
        $a_{n_j}\ne a_{n_{j+1}}, j=0, 1, 2, \cdots, k$ and 
        $\alpha_{n_0}\le \alpha_{n_1}\le \cdots\le \alpha_{n_k}$ and for all $\beta_{n_j},~ 0\le j\le n.$ Further let $M_{n_j}, 1\le j\le k $ with $M$ as defined in Theorem \ref{t1}.
    Then all the zeros lie in 
    \begin{eqnarray}
       |q|&\le& \frac{M}{|a_n|}\left\{ (\alpha_{n}-\alpha_{0}+|\alpha_{0}|)+2\sum_{j=0}^{k-1}|\beta_{n_j}|+|\beta_n|\right\}.
    \end{eqnarray}
\end{cor}

Setting $\beta_{n_j}=0,~~ j=0, 1, \ldots, n$ and taking $\alpha_{n_j}=a_{n_j}, ~~j=0, 1, \ldots, n$, we have the following result from Corollary \ref{co2}. 

\begin{cor}\label{cor3}
Let $p(q)=\sum_{\nu=0}^{n}q^{n_j}a_{n_j}$ such that no two adjacent coefficients $a_{n_j}'s$ are equal, i.e.
        $a_{n_j}\ne a_{n_{j+1}}, j=0, 1, 2, \cdots, k$ and 
        $\alpha_{n_0}\le \alpha_{n_1}\le \cdots\le \alpha_{n_k}$. Further let $M_{n_j}, 1\le j\le k $ with $M$ as defined in Theorem \ref{t1}.
Then all the zeros of $p(q)$ lie in 
    \begin{eqnarray}
       |q|&\le& \frac{M}{|a_n|}(\alpha_{n}-\alpha_{0}+|\alpha_{0}|).
    \end{eqnarray}
\end{cor}

\section{Conclusion}
Over the past ten years, the regular functions of the quaternionic variable have been introduced and thoroughly examined. Their diverse applications across a wide range of scientific disciplines have contributed significantly to their rapid development and have shown themselves to be a rich topic for investigation. We note that it became interesting to determine the regions containing some or all of the zeros of a regular polynomial of quaternionic variable following the study of the structure of zero sets and the Fundamental Theorem of Algebra for regular polynomials. The distribution of zeros for polynomials with quaternionic variables and quaternionic coefficients was not well covered in the literature.Here, we derive regions that contain all of the zeros of a regular polynomial of a quaternionic variable when the structure of the zero sets established in the recently developed theory of regular functions and polynomials of a quaternionic variable, as well as a maximum modulus theorem, restrict the real and imaginary parts of its coefficients.

%\footnotetext{Note: }


\end{document}